# ON THE SELBERG CLASS OF DIRICHLET SERIES: SMALL DEGREES

J. B. Conrey
A. Ghosh

## 1. Introduction.

In the study of Dirichlet series with arithmetic significance there has appeared (through the study of known examples) certain expectations, namely (i) if a functional equation and Euler product exists, then it is likely that a type of Riemann hypothesis will hold, (ii) that if in addition the function has a simple pole at the point $s = 1$, then it must be a product of the Riemann zeta-function and another Dirichlet series with similar properties, and (iii) that a type of converse theorem holds, namely that all such Dirichlet series can be obtained by considering Mellin transforms of automorphic forms associated with arithmetic groups. Guided by these ideas, consider the class $\mathcal{S}$ of Dirichlet series (introduced by Selberg [7]) : a Dirichlet series

$$F(s) = \sum_{n=1}^{\infty} \frac{a_n}{n^s}$$

is in $\mathcal{S}$ provided that it satisfies the following hypotheses:

(1) Analyticity: $(s-1)^m F(s)$ is an entire function of finite order for some non-negative integer $m$
(2) Ramanujan Hypothesis: $a_n \ll_\epsilon n^\epsilon$ for any fixed $\epsilon > 0$
(3) Functional equation: there must be a function $\gamma_F(s)$ of the form

$$\gamma_F(s) = \epsilon Q^s \prod_{i=1}^{k} \Gamma(w_i s + \mu_i)$$

where $|\epsilon| = 1$, $Q > 0$, $w_i > 0$, and $\Re \mu_i \geq 0$ such that

$$\Phi(s) = \gamma_F(s) F(s)$$

satisfies

$$\Phi(s) = \overline{\Phi}(1-s)$$

where $\overline{\Phi}(s) = \overline{\Phi(\overline{s})}$.

Research of the second author supported in part by the Alfred P. Sloan Foundation. Research of both authors supported in part by a grant from the NSF.

(4) Euler product: $a_1 = 1$, and

$$\log F(s) = \sum_{n=1}^{\infty} \frac{b_n}{n^s}$$

where $b_n = 0$ unless $n$ is a positive power of a prime and $b_n \ll n^\theta$ for some $\theta < 1/2$.

One is interested in classifying functions in this class and determining if properties (i) - (iii) hold. With this in mind it is clear that one is primarily interested in the primitive elements in this class namely the ones that cannot be written as a product of two or more non-trivial members of $\mathcal{S}$ (detailed definitions can be found in Section 2). For the problem of classification, one needs to define the degree $d_F$ of $F \in \mathcal{S}$ as

$$d_F = 2 \sum_{i=1}^{k} w_i,$$

allowing for the possibility $d_F = 0$ if the product of gamma functions is empty. Then $\mathcal{S}(d)$ will denote the subset of $\mathcal{S}$ consisting of all functions of degree $d$.

**Remarks**

(1) We do not consider generalised Dirichlet series here as the example $\zeta(\frac{s+1}{3})$ would be a candidate that would violate the Riemann hypothesis.
(2) The condition that there be at most one pole, and that at $s = 1$ is natural since if we expect condition (ii) to hold and if $F$ has a finite number of poles, then for each pole $F(s)$ would have the Riemann-zeta function, suitably shifted, as a factor. So the poles would lie on the line $\Re s = 1$ (otherwise $F$ would have zeros, corresponding to the shifted zeros of the Riemann zeta-function, not on the line $\Re s = \frac{1}{2}$). Thus we may write

$$F(s) = E(s) \prod_{j=1}^{R} \zeta(s + it_j)$$

for some real numbers $t_j$ and an entire function $E$. We then find nothing new by allowing these poles into our condition (1) and instead focus on functions with at most one pole (normalised to be at $s = 1$).
(3) In the functional equation, the restriction $\Re \mu_i \geq 0$ may be explained in the following way. Suppose there exists an arithmetic subgroup of $SL(2, R)$ together with a Maass cusp-form that corresponds to an exceptional eigenvalue and also assume that the Ramanujan-Petterson conjecture holds. Then the L-function that is associated with the Maass form has a functional equation with a $\mu_i$ satisfying $\Re \mu_i < 0$ but that violates the Riemann hypothesis. This suggests that a restriction of the type $\Re \mu_i \geq 0$ is appropriate.
(4) Condition (4) corresponds to the familiar notion of Euler product. In fact it is easy to verify that if $F \in \mathcal{S}$ then the coefficients $a_n$ of its Dirichlet series are multiplicative, i.e. $a_{mn} = a_m a_n$ if $m$ and $n$ are relatively prime.

Consequently, $F$ has an Euler product expansion

$$F(s) = \prod_p F_p(s) \qquad (\sigma > 1)$$

where the product is over primes and

$$F_p(s) = \sum_{k=0}^{\infty} \frac{a_{p^k}}{p^{ks}} \qquad (\sigma > 0).$$

From the viewpoint of automorphic L-functions, it is natural to put additional restrictions on $F_p(s)$ namely that for almost all primes $1/F_p(s)$ is a polynomial in $p^{-s}$ of degree independent of $p$. However, in this paper we will not need such a restriction.

(5) The condition $\theta < 1/2$ turns out to be surprisingly important. Note that $\theta = 1/2$ would violate the Riemann hypothesis as the example $(1-2^{1-s})\zeta(s)$ shows. It also plays a crucial part in Theorem 3.1 .

In Selberg [7] several conjectures were made concerning functions in $\mathcal{S}$. These are

(1) Regularity of distribution: There is an integer $n_F$ associated to each $F$ such that
$$\sum_{p \leq x} \frac{|a_p|^2}{p} = n_F \log \log x + \mathrm{O}(1)$$

(2) Orthonormality: If $F$ and $F'$ are distinct and primitive, then $n_F = 1$ and
$$\sum_{p \leq x} \frac{a_p a'_p}{p} = \mathrm{O}(1)$$

(3) $GL(1)$ twists: If $\chi$ is a primitive Dirichlet character and if $F$ and $F^\chi \in \mathcal{S}$ then
$$F = \prod_{i=1}^{k} F_i$$
where the $F_i$ are primitive implies that
$$F^\chi = \prod_{i=1}^{k} F_i^\chi$$
and the $F_i^\chi$ are also primitive elements of $\mathcal{S}$.

Selberg also conjectures that the Riemann hypothesis holds for this class of functions, i.e. if $F \in \mathcal{S}$, then all non-trivial zeros of $F(s)$ (see Sec. 2) are on the line $\sigma = 1/2$. However, we will not make use of this conjecture anywhere in this paper. In particular, in section 4 where we prove some consequences of Selberg's conjectures we do not use this one in our proofs.

Other fundamental questions also arise. For example, all instances of functions in $\mathcal{S}$ have degrees that are integers and so one would like to know if $\mathcal{S}(d)$ is empty when $d$ is not an integer. Similar questions arise as to the admissible values of $Q$ (see Sec. 3).

**Examples**

(1) Clearly $\zeta(s) \in \mathcal{S}(1)$ as is $L(s+i\alpha, \chi)$ for a primitive Dirichlet character $\chi$. Condition (1) prevents $\zeta(s+i\alpha)$ from being in the class if $\alpha \neq 0$.

(2) The Dedekind zeta-function of a number field of degree $d$ is in $\mathcal{S}(d)$ as are the abelian $L$- functions of primitive characters associated with such. Also, an Artin $L$-function for an irreducible representation of the Galois group of the number field is in $\mathcal{S}$ provided that it is holomorphic (i.e. that Artin's conjecture holds).

(3) An $L$-function associated with a holomorphic newform on a congruence subgroup of $SL_2(\mathbf{Z})$, once it is suitably normalized, is in $\mathcal{S}(2)$. An $L$-function associated with a non-holomorphic newform is presumably also in this class, but condition (2) has not been established for such. Also, if such a newform corresponds to an exceptional eigenvalue then condition (4) does not hold since $\Re\mu < 0$.

(4) The Rankin-Selberg convolution of any two holomorphic newforms is in $\mathcal{S}$. The symmetric square $L$-function associated with a holomorphic newform on the full modular group is in $\mathcal{S}(3)$.

(5) One can show that a large class of functions belong to $\mathcal{S}$, if one is willing to accept other well-known conjectures. For instance, Langland's conjectures allow L-functions associated with symmetric power representations to be in $\mathcal{S}$. Indeed, if $L(s)$ is associated with a cusp-form on the modular group (say) and
$$L(s) = \prod_p (1 - \frac{\alpha_p}{p^s})^{-1}(1 - \frac{\beta_p}{p^s})^{-1}$$
belongs in $\mathcal{S}$ and if $L_m(s)$ denotes the m-th symmetric power L-function
$$L_m(s) = \prod_p \prod_{j=0}^{m} (1 - \alpha_p^j \beta_p^{m-j}/p^s)^{-1},$$
then the Langlands' conjectures imply that $L_m(s)$ belongs to $\mathcal{S}$ for each $m > 0$, and they are entire. It was shown by Serre that this together with the non-vanishing on the line $\sigma = 1$ implies the Sato-Tate conjecture (on the distribution of the arguments of the $\alpha_p$'s). It can then be shown quite easily that the Sato-Tate conjecture implies that $n_{L_m} = 1$ for all m so that Selberg's conjectures would imply that all these L-functions are primitive.

It can be shown that the existence of a continuous distribution function for the arguments of the $\alpha_p$'s leads to the statement: the Sato-Tate conjecture is true if and only if $n_{L_m} = 1$ for all m. If one assumes only that the $n_{L_m}$ are integers, then one can show that there are at most two possible distributions, the Sato-Tate being one. There are similar statements for the Hasse-Weil L-functions (assuming the Taniyama-Weil conjectures) and Artin L-functions (assuming Artin's holomorphy conjecture).

In this paper, we discuss in some detail the structure of $\mathcal{S}(d)$ for $d \leq 2$, some consequences of Selberg's conjectures and the structure of a subclass of $\mathcal{S}(1)$ and $\mathcal{S}(2)$.

**2. Definitions.** We reserve $F$ and $G$ to denote members of $\mathcal{S}$.

For $F \in \mathcal{S}$ the $\epsilon$, $Q$, $w_i$, and $\mu_i$ are not well defined. For example, we could use the duplication formula

$$\Gamma(ws+\mu) = \Gamma(ws/2+\mu/2)\Gamma(ws/2+(1+\mu)/2)2^{1-ws}\pi^{-1/2}$$

to obtain a different $Q$, $w_i$, and $\mu_i$. Also, $\epsilon$ could be replaced by $-\epsilon$ without affecting the hypotheses for $F$ to be in $\mathcal{S}$. However, we shall show momentarily

**Theorem 2.1.** *If $\gamma^{(1)}$ and $\gamma^{(2)}$ are both admissible gamma factors for $F$, then $\gamma^{(1)}(s) = C\gamma^{(2)}(s)$ for some real constant $C$.*

Thus, $\gamma_F$ is well defined up to a constant. As this constant presents no particular problem for the moment we will use the notation $\gamma_F$ for one of the class of gamma factors. Later we will single out a particular choice. Now we define some notions.

**Trivial zeros:** The poles of $\gamma_F(s)$ are at the numbers $s = -(n+\mu_i)/w_i$, $1 \leq i \leq k$, $n = 0, 1, 2, \ldots$, which are all in $\sigma \leq 0$. From the functional equation it follows that $F$ has a zero at $s$ of order

$$-m_s + \sum_{\substack{i,n \\ n+\mu_i = w_i}} 1$$

where $m_s$ is the order of the pole of $F$ at $s$, i.e. $m_1 = m$ and $m_s = 0$ if $s \neq 1$. These zeros are the trivial zeros of $F$ and all other zeros are the non-trivial zeros. We do not exclude the possibility that $F$ has a trivial zero and a non-trivial zero at the same point.

**Degree:** We define the degree $d$ of a gamma factor for $F$ by the formula

$$d = 2\sum_{i=1}^{k} w_i.$$

We allow for the possibility that $d = 0$ as would happen if the product of gamma functions were empty. We will see later that the function

$$F(s) = 1$$

is the only degree 0 function in $\mathcal{S}$.

Note that if $\gamma^{(1)}(s)$ and $\gamma^{(2)}(s)$ are two admissible gamma factors for $F$ with degrees $d^{(1)}$ and $d^{(2)}$, then $d^{(1)} = d^{(2)}$. For, if we form the quotient of the two functional equations for $F$, say $h = \gamma^{(1)}/\gamma^{(2)}$, we obtain

$$h(s) = \overline{h}(1-s).$$

Now the left side is regular and non-zero in $\sigma > 0$ and the right side is regular and nonzero in $\sigma < 1$. Hence, $h(s)$ is an entire non-vanishing function. If the degrees were different, then $h$ would have zeros or poles. Thus, we may speak of the degree $d_F$ of $F$.

Arguing further, we give the proof of Theorem 2.1. We see that $h$ is of order 1 by Stirling's formula. Hence, it has the shape

$$e^{as+b}.$$

Then,

$$e^{as+b} = e^{\overline{a}(1-s)+\overline{b}}$$

so that $a = -\overline{a}$, i.e. $a = i\alpha$ where $\alpha$ is real. Letting $s = it$, we see again by Stirling's formula that

$$\left|\frac{h(it)}{h(-it)}\right| \to 1$$

as $t \to \infty$. Therefore, $\alpha = 0$. Thus, Theorem 2.1 follows.

**Primitive:** Next, we say that $F \in \mathcal{S}$ is primitive if $F = F_1 F_2$ with $F_1, F_2 \in \mathcal{S}$ implies that $F_1 = 1$ or $F_2 = 1$.

**Twists:** For a primitive Dirichlet character $\chi$ and an $F \in \mathcal{S}$ we define

$$F^\chi(s) := \sum_{n=1}^\infty \frac{a_n \chi(n)}{n^s}.$$

We end this section with some useful lemmas.

**Lemma 2.1.** *If $F(s) \in \mathcal{S}$, then $F_p(s)$ has no zeros in $\sigma > 1/2$ and $F(s)$ has no zeros in $\sigma > 1$.*

. For

$$\log F_p(s) = \sum_{k=1}^\infty \frac{b_{p^k}}{p^{ks}}$$

converges absolutely for $\sigma > 1/2$ by condition (4) and so has no zeros in $\sigma > 1/2$. Since the Euler product converges absolutely for $\sigma > 1$ there can be no zeros of $F$ in this half-plane.

**Lemma 2.2.** *Suppose that $F \in \mathcal{S}$ and $F(s)$ has a pole or zero at $s = 1 + i\alpha$ where $\alpha$ is real. Then the sums*

$$\sum_{p \leq x} \frac{a_p}{p^{1+i\alpha}}$$

*are unbounded as $x \to \infty$.*

*Proof.* We have

$$F(s) \sim C(s - (1+i\alpha))^m$$

as $s = \sigma + i\alpha \to 1^+ + i\alpha$ for some integer $m \neq 0$. Then

$$\log F(s) \sim m \log(\sigma - 1).$$

But by the bound on the $b_n$ in condition (4),

$$\log F(s) = \sum_{n=1}^{\infty} \frac{b_n}{n^s} = \sum_p \frac{a_p}{p^s} + \mathrm{O}\left(\sum_p \sum_{k=2}^{\infty} \frac{|b_{p^k}|}{p^{k\sigma}}\right)$$

$$= \sum_p \frac{a_p}{p^s} + \mathrm{O}(1)$$

where the sums over $p$ are for primes $p$. Thus,

$$\sum_p \frac{a_p}{p^s} \sim m \log(\sigma - 1)$$

as $\sigma \to 1^+$. Let

$$S(x) = \sum_{p \leq x} \frac{a_p}{p^{1+i\alpha}}.$$

Assume that $S(x)$ is bounded. Then

$$\sum_p \frac{a_p}{p^s} = \int_1^{\infty} x^{1-\sigma} dS(x)$$

$$= (\sigma - 1) \int_1^{\infty} S(x) x^{-\sigma}\, dx = \mathrm{O}(1)$$

which is a contradiction.

**3. Non-existence of functions with small degrees.** In this section we prove

**Theorem 3.1.** *If $F \in \mathcal{S}$ then $F = 1$ or $d_F \geq 1$.*

*Proof.* Suppose that $d_F < 1$. Then (using the notation $\int_{(a,b)}$ to mean the integral from $c - i\infty$ to $c + i\infty$ for any $c$ with $a < c < b$) we have

$$h(x) = \sum_{n=1}^{\infty} a_n e^{-2\pi n x}$$

$$= \frac{1}{2\pi i} \int_{(1,2)} (2\pi x)^{-s} \Gamma(s) F(s)\, ds$$

$$= \frac{P(\log x)}{x} + K(x)$$

where $P$ is a polynomial and

$$K(x) = \sum_{n=0}^{\infty} \frac{(-1)^n F(-n)(2\pi x)^n}{n!}$$

$$= \sum_{n=0}^{\infty} \frac{(-1)^n \gamma_F(n+1) F(n+1)(2\pi x)^n}{\gamma_F(-n) n!}$$

is an entire function of $x$ since

$$\frac{\gamma_F(n+1)}{\gamma_F(-n)n!} \ll n^{-(1-d)n}A^n$$

for some $A > 0$.

Thus we see that $h(x)$ is analytic in the plane slit along the negative real axis. But $h$ is periodic with period $i$ so in fact $h$ has no singularities on the real axis. Thus $h$ is entire. Let $H(z) = \sum_{n=1}^{\infty} a_n e(nz) = h(iz)$. Now if $y > 0$, then

$$a_n e^{-2\pi ny} = \int_0^1 H(x+iy)e(-nx)\,dx.$$

Both sides of this equation are entire functions of the complex variable $y$. We differentiate twice with respect to $y$ and set $y = 0$:

$$(2\pi n)^2 a_n = -\int_0^1 H''(x)e(-nx)\,dx \ll \int_0^1 |H''(x)|\,dx \ll 1.$$

Hence $a_n \ll n^{-2}$. But then the Dirichlet series for $F(s)$ is absolutely convergent for $\sigma > -1$. In particular, $F$ is uniformly bounded in $\sigma > -1/2$. But this is easily seen to be a contradiction if $d > 0$ as

$$F(1-s) = \frac{\Phi(s)}{\gamma_F(1-s)} \sim \frac{\gamma_F(s)}{\gamma_F(1-s)}$$

for $\sigma \geq \sigma_0 > 1$; by Stirling's formula

$$\left|\frac{\gamma_F(s)}{\gamma_F(1-s)}\right| \sim c(\sigma)t^{d(\sigma-1/2)}$$

for some $c(\sigma) > 0$ as $t \to \infty$.

In the case $d = 0$ we argue slightly differently. Since $H$ is entire we see that its Fourier series expansion is a power series in $e(z)$ which is entire so that it is convergent in the whole plane. This conclusion necessitates that the $a_n$ must be small. In fact, the $a_n$ will be so small that the Dirichlet series for $F$ will be absolutely convergent in the whole complex plane.

Then the functional equation for $F$ can be viewed as an identity between absolutely convergent Dirichlet series. We write the functional equation as

$$\sum_{n=1}^{\infty} a_n \left(\frac{Q^2}{n}\right)^s = \epsilon^* Q \sum_{n=1}^{\infty} \frac{\overline{a_n}}{n} n^s.$$

It follows that if $a_n \neq 0$ for some $n$, then $Q^2/n$ is an integer. Therefore, $Q^2$ is an integer, and $a_n \neq 0$ implies that $n \mid Q^2$. Thus, the Dirichlet series is really a Dirichlet polynomial. Now if $Q^2 = 1$, then $F = 1$. We assume that $q := Q^2 > 1$.

Note that $a_1 = 1$ implies that

$$|a_q| = Q.$$

Since $a_n$ is a multiplicative function, it follows that there is a prime $p$ and a positive integer $r$ such that $p^r \mid\mid q$ and
$$|a_{p^r}| \geq p^{r/2}.$$

Now the Euler factor corresponding to the prime $p$ is
$$F_p(s) = \sum_{j=0}^{r} \frac{a_{p^j}}{p^{js}}$$
and its logarithm is
$$\log F_p(s) = \sum_{j=0}^{\infty} \frac{b_{p^j}}{p^{js}}.$$

Replace $p^{-s}$ by $x$, $a_{p^j}$ by $A_j$, and $b_{p^j}$ by $B_j$. Then the Euler factor is
$$P(x) = \sum_{j=0}^{r} A_j x^j$$
and its logarithm is
$$\log P(x) = \sum_{j=0}^{\infty} B_j x^j.$$

Since $A_1 = 1$ we can factor $P$ as
$$P(x) = \prod_{i=1}^{r} (1 - R_i x).$$

Taking the logarithm of both sides here we obtain a formula for $B_j$:
$$B_j = -\sum_{i=0}^{r} \frac{R_i^j}{j}.$$

Now
$$\prod_{i=1}^{r} |R_i| = Q.$$

Therefore,
$$\max_{1 \leq i \leq r} |R_i| \geq p^{1/2}.$$

Then
$$|b_{p^j}|^{1/j} = |B_j|^{1/j} = \left| \sum_{i=1}^{r} \frac{R_i^j}{j} \right|^{\frac{1}{j}} \to \max_{1 \leq i \leq r} |R_i| \geq p^{1/2}.$$

This contradicts the existence of a $\theta < 1/2$ for which
$$b_n \ll n^\theta$$

and so completes the proof of the theorem.

This method of proof can be used to show that if $F \in \mathcal{S}$ and $d = 1$, then $Q \geq \pi^{-1/2}$. The modification is that if $d = 1$ and $Q < \pi^{-1/2}$ then $K$ is no longer entire but is regular in a circle about the origin of radius $1/2 + \eta(Q)$ where $\eta(Q)$ is a positive function of $Q$ which tends to 0 as $Q \to \pi^{-1/2}$. Then by periodicity we see that the power series $h(\exp(-2\pi x))$ is regular (apart from the negative real axis in a strip $\sigma > -\eta_1(Q)$ where $\eta_1(Q)$ is a positive function of $Q$ which tends to 0 as $Q \to \pi^{-1/2}$. Then the series for $h$ is convergent in this region which once again forces the $a_n$ to be too small to be coefficients of a function in $\mathcal{S}$.

We remark that Bochner [1] has a theorem which is relevant here. (See also Vigneras [10].) His result in our context is

**Theorem 3.2 (Bochner).** *Fix $Q > 0$, $w_i > 0$, $\mu_i \in \mathbf{C}$, for $1 \leq i \leq k$ with $\sum_{i=1}^{k} w_i = 1/2$. The number of linearly independent Dirichlet series $F(s)$ which satisfy a functional equation with a $\gamma_F$ satisfying these requirements is*

$$\leq 2q \prod_{i=1}^{k} w_i^{2w_i}$$

*where $q = \pi Q^2$.*

His proof involves Fuch's theorem on differential equations and Polya's gap theorem on singularities of power series.

We remark that as a consequence of Theorem 3.1 we have

**Corollary 3.3.** *Any function in $\mathcal{S}$ can be factored into a product of primitives.*

*Proof.* This assertion follows from the additivity of the degree function: $d_{FG} = d_F + d_G$ and the lower bound that $F \neq 1 \implies d_F \geq 1$.

As a second consequence we have

**Corollary 3.4.** *If $F \in \mathcal{S}$ and $d_F < 2$ then $F$ is primitive.*

**4. Consequences of Selberg's conjectures.** In this section we assume Selberg's conjectures and give some of the immediate consequences, many of which are mentioned in Selberg's paper.

**Proposition 4.1.** *If*

$$F = F_1^{e_1} F_2^{e_2} \ldots F_k^{e_k}$$

*where the $F_i$ are primitive, then*

$$n_F = e_1^2 + e_2^2 + \ldots e_k^2.$$

*Proof.* It is clear that

$$a_p(F) = \sum_{i=1}^{k} e_i a_p(F_i).$$

Then

$$|a_p(F)|^2 = \sum_{i=1}^{k} e_i^2 |a_p(F_i)|^2 + \sum_{i \neq j} a_p(F_i) \overline{a_p(F_j)}$$

whence the result follows by orthogonality.

**Proposition 4.2.** *If $n_F = 1$, then $F$ is primitive.*

*Proof.* This follows immediately from the first Proposition since if $F$ had a factorization as above then

$$n_F = 1 = \sum_{i=1}^{k} e_i^2$$

so that $k = 1$ and $e_1 = 1$.

**Proposition 4.3.** *If $F(s)$ has a pole of order $m$ at $s = 1$, then $\zeta(s)^m$ divides $F$.*

*Proof.* Clearly it suffices to prove this assertion in the case that $m = 1$ and $F$ is primitive. By Proposition 4.2, $\zeta$ is primitive since $n_\zeta = 1$. Take $F_1 = F$ and $F_2 = \zeta$ in the orthogonality conjecture. If $F \neq \zeta$, then that conjecture implies that $\sum_{p \leq x} \frac{a_p(F)}{p} = O(1)$. But this contradicts Lemma 2.2.

**Proposition 4.4.** *$\mathcal{S}$ has unique factorization.*

*Proof.* It suffices to show that if $F$ is primitive and $F \mid GG'$, then $F \mid G$ or $F \mid G'$. Assume neither holds. Suppose that $FF' = GG'$. We express both sides of this equation in terms of primitive functions as

$$F^f F_1^{f_1} \ldots F_k^{f_k} = G_1^{g_1} \ldots G_l^{g_l}$$

where $F$, $F_i$, and $G_i$ are distinct primitive functions. Multiply both sides by $F^r$ and compute $n_{()}$ for both sides:

$$(r + f)^2 + O(1) = r^2 + O(1).$$

We have a contradiction as $r \to \infty$.

**Proposition 4.5 (Dedekind's conjecture).** *If $\mathbf{K}$ is a finite extension of $\mathbf{Q}$ and $\zeta_{\mathbf{K}}$ is the Dedekind zeta function of $\mathbf{K}$, then*

$$L(s) = \zeta_{\mathbf{K}}(s)/\zeta(s)$$

*is entire.*

For $\zeta_{\mathbf{K}} \in \mathcal{S}$ and has a simple pole at $s = 1$. Hence, it will be divisible by $\zeta$ in $\mathcal{S}$.

We remark that R. Murty, in work to appear, has shown that Artin's conjecture about the holomorphy of Artin $L$-functions is also a consequence of Selberg's conjectures.

**Proposition 4.6.** *If $F \in \mathcal{S}$ then $F$ has no zeros on $\sigma = 1$.*

*Proof.* Clearly it suffices to prove the assertion for a primitive $F$. The assertion is true for $\zeta$, so we may assume that $F$ is entire. Then $F(s - i\alpha)$ is also a primitive member of $\mathcal{S}$. Applying the orthogonality relations to $F$ and $\zeta$ we see that $\sum_{p \leq x} \frac{a_p}{p^{1+i\alpha}} = O(1)$

**5. The class $\mathcal{S}^*(1)$.**

In all known examples of $F \in \mathcal{S}$ it is the case that one may find a $\gamma_F$ in which all $w_i = 1/2$. With this normalization, $Q$ is uniquely determined as are the $\mu_i$. Also, $\epsilon$ is ambiguous only as far as a factor of $\pm 1$. Thus, $\epsilon^2$ is uniquely determined. We are led to consider the possibly smaller class of functions $\mathcal{S}^*$ defined by the same axioms as $\mathcal{S}$ except that the functional equation has the form

$$\epsilon^2 Q^s \prod_{i=1}^{d} \Gamma(s/2 + \mu_i) F(s) = Q^{1-s} \prod_{i=1}^{d} \Gamma((1-s)/2 + \overline{\mu_i}) \overline{F}(1-s).$$

We note that to each member of $\mathcal{S}^*$ there is a unique 4-tuple

$$(d; q; \mu_1, \ldots, \mu_d; \epsilon^2)$$

where

$$q = \pi^d Q^2.$$

In practice $q$ is always a positive integer, though we do not make that assumption. We use the notation $\mathcal{S}^*(d)$ to denote all the elements of $\mathcal{S}^*$ with a given value of $d$. Thus, $\mathcal{S}^*(1)$ consists of all Dirichlet series in $\mathcal{S}^*$ with $d = 1$.

**Theorem 5.1.** *Suppose that $F \in \mathcal{S}^*(1)$. Then $q$ is a positive integer and there exists a primitive Dirichlet character $\chi$ mod $q$ and a real number $\alpha$ such that $F(s) = L(s + i\alpha, \chi)$.*

Our proof is similar to Siegel's proof (see [8] or [9]) of Hamburger's theorem, though it is formulated somewhat differently. Also, Gerardin and Li [3] have given an adelic proof of a similar theorem. Note the important role played by the bound $\theta < 1/2$ of axiom (4) of the definition of the Selberg class. If not for that condition the Dirichlet series

$$\sum_{n \text{ odd}} n^{-s} - 2 \sum_{n \equiv 2 \mod 4} n^{-s} = (1 - 2^{1-s}) \prod_{p \geq 3} (1 - p^{-s})^{-1}$$

would be an element of $\mathcal{S}$ even though this function does not satisfy the Riemann Hypothesis. Thus, (4) features into our proof in a significant way.

*Proof.* We assume that $F$ has a functional equation

$$\epsilon Q^s \Gamma\left(\frac{s}{2} + \mu\right) F(s) = \overline{\epsilon} Q^{1-s} \Gamma\left(\frac{1-s}{2} + \overline{\mu}\right) \overline{F}(1-s).$$

If $\mu$ is not real then we consider

$$G(s) = F(s + \overline{\mu} - \mu).$$

Then $G$ satisfies the functional equation

$$\epsilon_1 Q^s \Gamma\left(\frac{s}{2} + \frac{\mu + \overline{\mu}}{2}\right) G(s) = \overline{\epsilon_1} Q^{1-s} \Gamma\left(\frac{1-s}{2} + \frac{\mu + \overline{\mu}}{2}\right) \overline{G}(1-s)$$

where
$$\epsilon_1 = \epsilon Q^{\overline{\mu}-\mu}.$$

Also, $G(s)$ has a pole possibly at $s = 1 + \mu - \overline{\mu}$. Thus, without loss of generality we may assume that $\mu$ is real and that the pole of $F$ (if it exists) is at $s = 1 + i\alpha$ for some real $\alpha$.

Let $q = \pi Q^2$. Then we consider
$$f(x) = \sum_{n=1}^{\infty} a_n \left(\frac{2\pi n x}{q}\right)^{2\mu} e^{-2\pi n x/q} = \frac{1}{2\pi i} \int_{(1,2)} (2\pi x/q)^{-s} \Gamma(s + 2\mu) F(s) \, ds.$$

We move the line of integration to the left of 1 and get a residue from the pole at $s = 1 + i\alpha$. We make the change of variables $s \to 1 - s$ in the integral and apply the functional equation. Thus,

$$f(x) = \frac{P_1(\log x)}{x^{1+i\alpha}}$$
$$+ \frac{1}{2\pi i} \int_{(0,1)} (2\pi x/q)^{s-1} \Gamma(1 - s + 2\mu) \overline{\epsilon}^2 Q^{2s-1} \frac{\Gamma(s/2 + \mu)}{\Gamma((1-s)/2 + \mu)} \overline{F}(s) \, ds.$$

We make the change of variables $s \to 2s$ and use the duplication formula on the first gamma function in the integral. Then the integral is

$$2Q\overline{\epsilon}^2 \int_{(0,1/2)} (2x)^{2s-1} \frac{\Gamma(1 + 2\mu - 2s)\Gamma(s + \mu)}{\Gamma(1/2 + \mu - s)} \overline{F}(2s) \, ds$$

$$= \frac{Q\overline{\epsilon}^2 2^{2\mu}}{\pi^{1/2} x} \int_{(0,1/2)} \frac{\Gamma(1/2 + \mu - s)\Gamma(1 + \mu - s)\Gamma(s + \mu)}{\Gamma(1/2 + \mu - s)} x^{2s} \overline{F}(2s)$$

$$= \frac{C_q}{x} \int_{(0,1/2)} \Gamma(1 + \mu - s)\Gamma(s + \mu) x^{2s} \overline{F}(2s)$$

where
$$C_q = \frac{Q\overline{\epsilon}^2 2^{2\mu}}{\pi^{1/2}}.$$

Now we move the line of integration to the right of $1/2$ crossing the pole at $s = 1/2 - i\alpha/2$. Then we are in a region where the Dirichlet series $F(2s)$ converges absolutely. We expand $F(2s)$ into its Dirichlet series and integrate term-by-term. We have

$$f(x) = \frac{P_1(\log x)}{x^{1+i\alpha}} + x^{-i\alpha} P_0(\log x) + C_q \Gamma(1 + 2\mu) x^{-1-2\mu} \sum_{n=1}^{\infty} \frac{\overline{a_n} n^{2\mu}}{(1 + n^2/x^2)^{1+2\mu}}.$$

The latter formula comes from the integral formula

$$\frac{1}{2\pi i} \int_{(-b,a)} \Gamma(s+b)\Gamma(a-s) x^{-s} \, ds = \frac{x^b \Gamma(a+b)}{(1+x)^{a+b}}.$$

We divide our equation by $x^{2\mu}$ and obtain

$$Per(x) = L(x) + Cx \sum_{n=1}^{\infty} \frac{\overline{a_n} n^{2\mu}}{(x^2 + n^2)^{1+2\mu}}$$

where $Per(x)$ is a function that is regular in $\Re x > 0$ and is periodic with period $iq$. Also, $L(x)$ is a function which is regular in the whole plane with the negative real axis removed and $C$ is independent of $x$. Now as $x = \eta + in \to 0^+ + in$, the right side is asymptotic to

$$C'\overline{a_n}$$

for some $C'$ which is independent of $n$. By periodicity, we have the same asymptotics as $x = \eta + i(n+q) \to 0^+ + i(n+q)$. Therefore, $q$ must be an integer and $a_n = a_{n+q}$.

Now we use the fact that $a_n$ is multiplicative. It is not difficult to show that if $a_n$ is a multiplicative function which is periodic mod $q$, then there exists a Dirichlet character $\chi_1$ mod $q$ for which $a_n = \chi_1(n)$ for all $n$ for which $(n, q) = 1$. It suffices to show that $a_n$ is completely multiplicative on such $n$. So suppose that $(mn, q) = 1$. Let $r$ be an integer for which $(m + rq, n) = 1$. Then

$$a_m a_n = a_{m+rq} a_n = a_{(m+rq)n} = a_{mn}$$

whence $a_n$ is completely multiplicative on those $n$ with $(n, q) = 1$. Now let $\chi$ be the primitive character which induces $\chi_1$. Suppose that $\chi$ is a character with modulus $q_1$. Then $L(s, \chi)$ satisfies a functional equation

$$\epsilon_\chi \left(\sqrt{\frac{q_1}{\pi}}\right)^s \Gamma((s+a)/2) L(s, \chi) = \Phi(s) = \overline{\Phi}(1-s).$$

If we form the quotient of the functional equations for $F(s)$ and $L(s, \chi)$ we obtain an equation

$$\phi(s) := \frac{\Gamma(s/2 + \mu)}{\Gamma(s/2 + a/2)} E(s) = CQ_1^s \frac{\Gamma((1-s)/2 + \mu)}{\Gamma((1-s)/2 + a/2)} \overline{E}(1-s)$$

for some $C$ which is independent of $s$, and where

$$E(s) = F(s)/L(s, \chi) = \prod_{p | q_1} F_p(s).$$

Also, $a$ is either 0 or 1. We show that such an equation can only hold if $\mu = a/2$ and $E(s) = 1$.

To do this we first show that all zeros and poles of $\phi(s)$ are in $\sigma \leq 1/2$. For the quotient of the gamma functions, all zeros and poles are in $\sigma \leq 0$. That $E(s)$ has no zeros in $\sigma > 1/2$ follows from the fact that its Euler product involves only a finite number of factors, each of which has no zeros in $\sigma > 1/2$ by Lemma 1.2.

By the symmetry of the functional equation it follows that all zeros and poles of $\phi(s)$ are on $\sigma = 1/2$. But then the quotient of the gamma functions is entire with no zeros. Hence, the gamma functions are the same, i.e. $\mu = a/2$. Then we are left with a degree 0 functional equation, which we've seen in Theorem 3.1 implies that the $b_n$ are too large. Hence, $\phi(s) = 1$, and $F(s) = \zeta(s)$ or $F(s) = L(s, \chi)$ for a primitive character $\chi$.

**6. Functional equations from $\mathcal{S}^*(2)$.** We now prove a general "converse" theorem about Dirichlet series which have GL(2) type functional equations. This theorem may be regarded as a generalization of the basic theorems of Hecke [5] and Maass [6], and contains them as special cases. We set up some notation for this section only. Let

$$F(s) = \sum_{n=1}^{\infty} \frac{a_n}{n^s}$$

be absolutely convergent for $\sigma > 1$. Let

$$\Phi(s) = (\sqrt{q}/\pi)^s \Gamma\left(\frac{\alpha + \beta + 1/2 + s}{2}\right) \Gamma\left(\frac{\alpha - \beta + 1/2 + s}{2}\right) F(s).$$

We assume that $\alpha$ is real as this presents no loss of generality. Let $J$ and $K$ denote the usual Bessel functions, and define

$$H_\alpha(x) = x^{1/2} J_\alpha(x).$$

For $x$ and $y$ positive let

$$f(x, y) = y^{1/2} \sum_{n=1}^{\infty} a_n H_\alpha(2\pi n x/\sqrt{q}) K_\beta(2\pi n y/\sqrt{q}).$$

**Theorem 6.1.** $\Phi(s)$ is entire and satisfies

$$\Phi(s) = \overline{\Phi}(1 - s)$$

if and only if

$$f\left(re^{i\theta}\right) = \overline{f\left(r^{-1} e^{i\theta}\right)}$$

where $0 < \theta < \pi/2$, $r > 0$, and

$$re^{i\theta} = x + iy.$$

We remark that $H_{-1/2}(x) = (2/\pi)^{1/2} \cos x$ and $H_{1/2}(x) = (2/\pi)^{1/2} \sin x$. In fact, $\alpha = -1/2$ in the above theorem corresponds to the case of even Maass forms while $\alpha = 1/2$ corresponds to the situation of odd Maass forms. The work of Epstein, Hafner, and Sarnak [2] already contains this case. Our theorem may be regarded as a generalization of their result. When $\beta = 1/2$ we are in the situation of holomorphic cusp forms. To see this, we first observe that $y^{1/2} K_{1/2}(y) = (2\pi)^{-1/2} e^{-y}$. Next, we note that if $\alpha$ is half an odd integer, say $\alpha = (k-1)/2$ where $k$ is an even integer, then $H_\alpha(x)$ has a zero at $x = 0$ of order $k/2$. Then differentiating $k/2$ times with respect to $x$ and setting $x = 0$ we obtain the usual representation for the modular form. On the other hand, when we differentiate $k/2$ times the relation

$$f(x, y) = f\left(\frac{x}{x^2 + y^2}, \frac{y}{x^2 + y^2}\right)$$

with respect to $x$ and set $x = 0$ and use the fact that the derivatives lower than $k/2$ vanish at $x = 0$ we obtain that for some constant $c$,

$$g(y) = c \left(\frac{\partial}{\partial x}\right)^{k/2} f(x,y)|_{x=0} = \sum_{n=1}^{\infty} a_n n^{(k-1)/2} e^{-2\pi n y}$$

satisfies
$$g(1/y) = y^k g(y)$$

which is the desired relation. As an example, consider the case of the $\Delta$-function generated by the Ramanujan coefficients $\tau(n)$:

$$\Delta(z) = \sum_{n=1}^{\infty} \tau(n) e(nz).$$

This corresponds to the situation $\beta = 1/2$, $\alpha = 11/2$. We have

$$f(x,y) = x^{1/2} \sum_{n=1}^{\infty} \frac{\tau(n)}{n^{11/2}} J_{11/2}(2\pi n x) e^{-2\pi n y}$$

and
$$F(s) = \sum_{n=1}^{\infty} \frac{\tau(n)}{n^{11/2}} n^{-s}$$

and
$$f(x,y) = f\left(\frac{x}{x^2+y^2}, \frac{y}{x^2+y^2}\right)$$

if and only if (after use of the duplication formula)

$$(2\pi)^{-s} \Gamma(s + 11/2) F(s) = \Phi(s) = \Phi(1-s).$$

The latter is, of course, the well known functional equation (after a shift by $11/2$) for the $L$-function associated with $\tau$. The former is less recognizable. One calculates that

$$(\pi/2)^{1/2} x^{1/2} J_{11/2}(x) = \cos x \left(-1 + \frac{105}{x^2} - \frac{945}{x^4}\right) + \sin x \left(\frac{15}{x} - \frac{420}{x^3} + \frac{945}{x^5}\right)$$

and that this function has a zero of order 6 at $x = 0$. After differentiating 6 times with respect to $x$ and setting $x = 0$ one obtains the usual transformation formula for the $\Delta$ function

$$\Delta(iy) = y^{-12} \Delta(1/iy).$$

We remark that $f(x,y)$ above satisfies the differential equation

$$f_{xx}(x,y) + f_{yy}(x,y) = \left(\frac{\alpha^2 - 1/4}{x^2} + \frac{\beta^2 - 1/4}{y^2}\right) f(x,y).$$

Also, the above theorem remains valid if we replace the Dirichlet series by a generalized Dirichlet series

$$F(s) = \sum_{n=1}^{\infty} a_n \lambda_n^{-s}$$

and replace the Bessel series by

$$f(x,y) = (xy)^{1/2} \sum_{n=1}^{\infty} a_n J_\alpha(2\pi\lambda_n x/\sqrt{q}) K_\beta(2\pi\lambda_n y/\sqrt{q})$$

for a fairly general sequence $\lambda_n$.

We begin with the following lemma about the Mellin transform of a product of Bessel functions:

**Lemma 6.2.** *Suppose that $a$ and $b$ are real and positive. Then for $\Re s > 0$ we have*

$$\int_0^\infty J_\alpha(au) K_\beta(bu) u^s \frac{du}{u} = \left(\frac{a}{b}\right)^\alpha b^{-s} 2^{s-2}$$

$$\cdot \frac{\Gamma\left(\frac{s+\alpha+\beta}{2}\right) \Gamma\left(\frac{s+\alpha-\beta}{2}\right)}{\Gamma(\alpha+1)} {}_2F_1\left(\frac{s+\alpha+\beta}{2}, \frac{s+\alpha-\beta}{2}, \alpha+1, -\frac{a^2}{b^2}\right).$$

*Proof of Lemma 6.2.* This lemma follows easily from the transforms

$$\int_0^\infty J_\alpha(x) x^s \frac{dx}{x} = 2^{s-1} \Gamma\left(\frac{s+\alpha}{2}\right) / \Gamma\left(\frac{\alpha+2-s}{2}\right)$$

and

$$\int_0^\infty K_\beta(x) x^s \frac{dx}{x} = 2^{s-2} \Gamma\left(\frac{s+\beta}{2}\right) \Gamma\left(\frac{s-\beta}{2}\right)$$

and

$$\frac{1}{2\pi i} \int_{(c)} \frac{\Gamma(s)\Gamma(l-s)\Gamma(m-s)}{\Gamma(n-s)} x^{-s} ds = \frac{\Gamma(l)\Gamma(m)}{\Gamma(n)} {}_2F_1(l,m,n,-x)$$

for $\Re s > 0$ by using the convolution property of Mellin transforms:

$$\int_0^\infty f(x) g(x) x^s \frac{dx}{x} = \frac{1}{2\pi i} \int_{(c)} F(z) G(s-z) dz$$

where $f$ and $F$ are a Mellin transform pair as are $g$ and $G$.

**Lemma 6.3.** *Let*

$$T(s) = (\sin\theta)^{-s} {}_2F_1\left(\frac{s+\alpha+\beta+1/2}{2}, \frac{s+\alpha-\beta+1/2}{2}, \alpha+1, -(\cot\theta)^2\right)$$

*for $0 < \theta < \pi/2$. Then $T(s) = T(1-s)$.*

*Proof of Lemma 6.3.* This fact follows easily from the transformation formula

$${}_2F_1(a,b,c,x) = (1-x)^{c-a-b} {}_2F_1(c-a, c-b, c, x)$$

and from the fact that $_2F_1(a,b,c,x) = {_2F_1}(b,a,c,x)$.

*Proof of Theorem 6.1.* We consider the Mellin transform

$$M(s) = \int_0^\infty f\left(re^{i\theta}\right) r^{s-1/2} \frac{dr}{r}.$$

Replacing $r$ by $1/r$ in the integral we see that

$$f\left(re^{i\theta}\right) = \overline{f\left(r^{-1}e^{i\theta}\right)}$$

implies that $M(s) = \overline{M}(1-s)$. By the definition of $f$ we find that

$$M(s) = (2\pi \sin\theta \cos\theta)^{1/2} \sum_{n=1}^\infty a_n n^{1/2}$$
$$\cdot \int_0^\infty J_\alpha((2\pi nr \cos\theta)/\sqrt{q}) K_\beta((2\pi nr \cos\theta)/\sqrt{q}) r^{s+1/2} \frac{dr}{r}.$$

We evaluate the integral with the help of Lemma 6.2 and find that

$$M(s) = 2^{-3/2} (\cos\theta)^{1/2} (\cot\theta)^\alpha \Gamma(1+\alpha)^{-1} T(s)\Phi(s).$$

Since $\alpha$ is real it follows from Lemma 6.3 that $\Phi(s) = \overline{\Phi}(1-s)$ if and only if $M(s) = \overline{M}(1-s)$. But by Mellin inversion

$$f\left(re^{i\theta}\right) = \frac{1}{2\pi i} \int_{(c)} M(s) r^{1/2-s} ds$$

so that replacing $s$ by $1-s$ we find that $M(s) = \overline{M}(1-s)$ implies that $f\left(re^{i\theta}\right) = \overline{f\left(r^{-1}e^{i\theta}\right)}$. This completes the proof.

**Acknowledgments.** The first author would like to thank the Rutgers Mathematics Department for providing an excellent place to work during his sabbatical leave when this manuscript was prepared. He especially thanks Bill Duke and Henryk Iwaniec for their hospitality during his stay and the many stimulating conversations regarding this work.

Both authors also thank the Institute for Advanced Study for providing a stimulating work environment during their sojourn to New Jersey.

## References

1. S. Bochner, *On Riemann's functional equation with multiple gamma factors*, Annals of Math. **67** (1958), 29–41.
2. C. Epstein, J. Hafner, and P. Sarnak, *Zeros of L-functions attached to Maass forms*, Math. Zeitschrift **190** (1985), 113–128.
3. P. Gerardin and W. Li, *Functional equations and periodic sequences*, Théorie des Nombres (Quebec, PQ, 1987), deGruyter, Berlin New York, 1989.
4. M. I. Gurevich, *Determining L-series from their functional equations*, Math. USSR Sbornik **14** (1971), 537–553.

5. E. Hecke, *Über die Bestimmung Dirichletscher Reihen durch ihre Funktionalgleichung*, Mathematische Annalen **112** (1936), 664–699.
6. H. Maass, *Über eine neue Art von nichtanalytischen automorphen Funktionen und die Bestimmung Dirichletscher Reihen durch Funktionalgleichungen*, Mathematische Annalen **121** (1949), 141–183.
7. A. Selberg, *Old and new conjectures and results about a class of Dirichlet series*, Collected Papers, vol. 2, Springer-Verlag, Berlin Heidelberg New York, 1991.
8. C. L. Siegel, *Bemerkungen zu einem Satz von Hamburger über die Funktionalgleichung der Riemannsche Zetafunktion*, Mathematische Annalen **86** (1922), 276–279.
9. E. C. Titchmarsh, *The Theory of the Riemann Zeta-function*, Second Edition, Clarendon Press, Oxford, 1986.
10. M.-F. Vigneras, *Facteurs gamma et équations fonctionnelles*, Springer-Verlag Lecture Notes 627, 1977, pp. 79–103.

Department of Mathematics, Oklahoma State University, Stillwater, Oklahoma 74078